\documentclass[10pt]{article}
\usepackage{amsmath,amsthm,amsfonts}
\usepackage{hyperref}
\usepackage[usenames]{color}

\newtheorem{theorem}{Theorem}[section]
\newtheorem{lemma}[theorem]{Lemma}
\theoremstyle{definition}
\newtheorem{algorithm}[theorem]{Algorithm}
\newtheorem{corollary}[theorem]{Corollary}
\newtheorem{conjecture}[theorem]{Conjecture}
\newtheorem{definition}[theorem]{Definition}
\newtheorem{example}[theorem]{Example}
\newtheorem{problem}[theorem]{Problem}
\newtheorem{remark}[theorem]{Remark}

\def\la{{\overline{a}}}

\newcommand{\gp}[1]{{\langle #1 \rangle}}

\newcounter{comcount}

\title{Random subgroups and analysis \\ of the length-based and quotient attacks }

\author{Alexei G. Myasnikov and Alexander Ushakov}

\date{July  09, 2007}

\begin{document}
\maketitle

\begin{abstract}
In this paper we discuss generic properties of "random subgroups" of
a given group $G$. It turns out that in many groups $G$ (even in
most exotic of them) the  random subgroups have a simple algebraic
structure and they "sit" inside $G$ in a very particular way. This
gives a strong mathematical foundation for cryptanalysis of several
group-based cryptosystems and  indicates on how to chose "strong
keys". To illustrate our technique we analyze the
Anshel-Anshel-Goldfeld (AAG) cryptosystem and give a mathematical
explanation of recent success of some heuristic length-based attacks
on it. Furthermore, we design and analyze a new type of attacks,
which we term the quotient attacks.  Mathematical methods we develop
here also indicate how one can try to choose "parameters" in AAG to
foil the attacks.
\end{abstract}

\tableofcontents

\section{Introduction}

Most of the modern cryptosystems use algebraic structures as their
platforms such as  rings, groups, lattices, etc.  Typically,
cryptographic protocols involve a random choice of various algebraic
objects related to the platforms: elements, subgroups, or
homomorphisms. One of the key points to use randomness is to foil
various statistical  attacks, or attacks  which could use some
specific properties of objects if they are not chosen randomly. The
main goal of this paper is to show that randomly chosen objects
quite often have very particular properties, which allow some
"unexpected" attacks. We argue that knowledge of basic properties of
the random objects must be a part of any serious cryptanalysis and
it has to be one of the principal tools in choosing  good keys.

In the paper \cite{MU1} we studied asymptotic properties of words
representing the trivial element in a given finitely presented group
$G$. It turned out that a randomly chosen trivial word in $G$ has a
"hyperbolic" van Kampen diagram, even if the group $G$ itself is not
hyperbolic. This allows one to design a correct (no errors) search
decision algorithm which gives the  answer in polynomial time on a
generic subset (i.e., on "most" elements) of the Word Search Problem
in $G$. A similar result for the Conjugacy Search Problem in
finitely presented groups has been proven in \cite{MU2}. These
results show that the group-based cryptosystems whose security is
based on the word or  conjugacy search problems are subject to
effective attacks, unless the keys are chosen in the complements of
the corresponding generic sets. Rigorous proofs of results of
\cite{MU1} and \cite{MU2} are available in \cite{U}.

In this paper we study  asymptotic properties of finitely generated
subgroup of groups. We start by introducing  a methodology to deal
with  asymptotic properties of subgroups in a given finitely
generated group, then we describe  two  such properties, and finally
we show how one can use them in cryptanalysis of group based
cryptosystems.

Then  we dwell on the role of asymptotically dominant properties of
subgroups in modern cryptanalysis. We mostly focus on one particular
example - the AAG cryptosystem \cite{AAG}, however, it seems
plausible that a similar analysis applies to some other
cryptosystems. One of our main goals here is to give mathematical
reasons why the so-called Length Based Attacks give surprisingly
good results in breaking AAG. Another goal is to introduce and
analyze a  new attack that we coiled  a {\em quotient attack}. We
also want to emphasize that we believe that this "asymptotic
cryptanalysis" provides a good method to choose strong keys (groups,
subgroups, and elements) for AAG scheme (with different groups as
the platforms) that may  prevent some of the known attacks,
including the ones discussed here.

 The main focus is on  security of the
Anshel-Anshel-Goldfeld (AAG) public key exchange scheme \cite{AAG}
and cryptanalysis of the Length Based Attack (LBA). This attack
appeared first in  the paper \cite{HT} by J. Hughes and A.
Tannenbaum, and later was further developed in a joint paper
\cite{GKTTV} by Garber, Kaplan, Teicher, Tsaban, and Vishne.
Recently, the most successful variation of this attack for braid
groups was developed in \cite{MU3}. Notice that Ruinsky, Shamir, and
Tsaban used LBA in attacking some other algorithmic problems in
groups \cite{ST}. Our goal is to give mathematical reasons  why
Length Based Attacks, which are, in their basic forms, very simple
algorithms,  give surprisingly good results in breaking AAG scheme.
It seems plausible that a similar analysis applies to some other
cryptosystems.  We hope that this cryptanalysis provides also a good
method to choose strong keys (groups, subgroups, and elements) for
various realizations of AAG schemes that would  prevent some of the
known attacks.

The basic idea of LBA is very simple, one solves the Simultaneous
Conjugacy Search Problem relative to a subgroup (SCSP*) (with a
constraint that the solutions are in a given subgroup) precisely in
the same way as this would be done in a free group.  Astonishingly,
experiments show that this strategy works well in groups which are
far from being free, for instance, in braid groups. We claim that
the primary reason for such phenomenon is that asymptotically
finitely generated subgroups in many groups are free. Namely, in
many  groups a randomly chosen tuple of elements with overwhelming
probability freely generates a free subgroup (groups with Free Basis
Property). This allows one to analyze the generic complexity of LBA,
SCSP*, and some other related algorithmic problems. Moreover, we
argue that LBA implicitly relies on fast computing of the geodesic
length of elements in finitely generated subgroups of the platform
group $G$, or some good approximations of that length.   In fact,
most of LBA strategies tacitly assume that the geodesic length of
elements in $G$ is a good approximation of the geodesic length of
the same elements in a subgroup. On the first glance this is a
provably wrong assumption,
 it is known that even in a braid group $B_n$, $n \geq 3$, there are infinitely many
 subgroups whose distortion function (that measures the geodesic length
 in a subgroup relative to the one in $G$) is  not bounded by any recursive function.
  We show, nevertheless, that, again, in many groups the distortion of randomly chosen finitely
  generated subgroups is at most linear. Our prime objective is the
  braid group $B_n, n \geq 3$. Unfortunately, the scope of this
  paper does not allow a thorough investigation of asymptotic
  properties of subgroups of $B_n$. However, we prove the main
  results for the pure braid groups $PB_n$, which are subgroups of
  finite index in the ambient braid groups. We conjecture that the
  results hold in the groups $B_n$ as well, and hope to fill in this
  gap elsewhere in the future. In fact, our results hold for all
  finitely generated groups $G$ that have non-abelian free
  quotients.

While studying the length based attacks we realized that there
exists a new powerful type of attacks on AAG cryptosystems - the
quotient attacks (QA). These  attacks are just fast generic
algorithms to solve various  search problems in groups, such as the
Membership Search  Problem (MSP) and  SCSP*. The main idea behind QA
is that to solve a computational  problem in a group $G$ it
suffices, on most inputs, to solve it in a suitable quotient $G/N$,
provided $G/N$ has a fast decision algorithm for the problem.
Robustness of such an algorithm relies on the following property of
the quotient $G/N$: a randomly chosen finitely generated subgroup of
$G$ has trivial intersection with the kernel $N$. In particular,
this is the case, if $G/N$ is a free non-abelian group. Notice, that
a similar idea was already exploited in \cite{KMSS1}, but there the
answer was given only for inputs in "No" part of a given decision
problem, which, obviously, does not apply to search problems at all.
The strength of our approach comes from the extra requirement that
$G/N$ has the free basis property.

  More generally, our main goal concerns with the  methods on how to use asymptotic
  algebra and generic case complexity in cryptanalysis of group based
  cryptosystems. All asymptotic results on subgroups, that are used here,  are based on
  the notion of an asymptotic density with respect to the standard
  distributions  on generating sets of the subgroups.
   Essentially, this notion  appeared first in the form of zero-one laws in probability
  theory and combinatorics. It became extremely popular after seminal works of Erdos,
  that shaped up the so-called The Probabilistic Method (see, for
  example, \cite{AS}).
  In infinite group theory it is due mostly to the famous Gromov's result on hyperbolicity
  of random finitely generated groups (see \cite{Ol} for a complete
  proof). Generic complexity of algorithmic problems appeared first in the
  papers \cite{KMSS1,KMSS2, BMR,BMS}. We refer the reader to a comprehensive survey \cite{GMMU}
   on generic complexity of algorithms.   Some recent relevant results on generic complexity of
   search problems in groups (which are of the main interest in
   cryptography) can be found in  \cite{U}.

This paper is intended to both algebraists and cryptographers. We
believe, that AAG cryptosystem,  despite being heavily battered by
several attacks, is very much alive still. It simply did not get a
fair chance to survive  because of insufficient group theoretic
research it required. It is still quite plausible   that there are
platform groups $G$ and methods to chose strong keys for AAG which
would foil all known attacks.  To find such a group $G$ is an
interesting algebraic  problem. On the other hand, our method of
analyzing generic complexity of computational security assumptions
of AAG, which is based on the asymptotic behavior of subgroups in a
given group, creates a bridge between  asymptotic algebra and
cryptanalysis. This could be applicable to other cryptosystems which
rely on a random choice of algebraic objects: subgroups, elements,
or homomorphisms.

\section{Asymptotically dominant properties}
\label{se:asym-properties}

In this section we develop some tools to study asymptotic properties
of subgroups of groups. Throughout this section by  $G$ we denote a
group with a finite generating set $X$.

\subsection{A brief description}
\label{subsec:random-subgroups}

Asymptotic properties of subgroups, a priori,   depend on a given
probability distribution on these subgroups. In general, there are
several natural procedures to generate subgroups in a given group.
However, there is no a unique universal distribution of this kind.
We refer to \cite{MO} for a discussion on different approaches to
random subgroup generation.

Our basic principle here is that in applications one has to consider
the particular  distribution that comes from a particular random
generator of subgroups used in the given application, say a
cryptographic protocol. As soon as the distribution is fixed one can
approach asymptotic properties of subgroups via asymptotic densities
 with respect to a fixed stratification of the set of subgroups
 (which usually comes alone with the generating procedure).  We
 briefly discuss these ideas below and refer to \cite{BMS,BMR,KMSS1, KMSS2},
 and to a recent survey
\cite{GMMU}, for a  thorough exposition. In Section
\ref{subsec:random-sub-gen} we adjust these general ideas to a
particular way to generate subgroups which is used  in cryptography.

 Recall, that  $G$ is  a group
generated by a finite set $X$. The first step is to choose and fix a
particular way to describe  finitely generated subgroups $H$ of $G$.
For example, a description $\delta$ of $H$ could be  a tuple of
words $(u_1, \ldots,u_k)$ in the alphabet $X^{\pm 1} = X \cup
X^{-1}$ representing a set of generators of $H$, or a set  of words
$\{u_1, \ldots,u_k\}$  that generates $H$, or a folded finite graph
that accepts the subgroup generated by the generators $\{u_1,
\ldots,u_k\}$ of $H$ in the ambient free group $F(X)$ (see
\cite{KM}), etc. In general, the descriptions above, by no means are
unique for a given subgroup $H$, in fact, we listed them here in the
decreasing degree of repetition.

When the way to describe subgroups in $G$ is fixed one can consider
the set $\Delta$ of all such descriptions of all finitely generated
subgroups of $G$. The next step is to define a {\em size}
$s(\delta)$ of a given description $\delta \in \Delta$, i.e., a
function
    $$s: \Delta \rightarrow  \mathbb{N}$$
in such a way that the set (the ball of radius $n$)
    $$B_n = \{\delta \in \Delta \mid s(\delta) \leq n\}$$
 is finite. This gives a {\em stratification} of the set $\Delta$ into a
union of finite balls:
 \begin{equation}
 \label{eq:stratification}
  \Delta = \cup_{n = 1}^{\infty} B_n.
   \end{equation}

  Let  $\mu_n$ be a given probabilistic measure  on $B_n$ (it could be the
  measure induced  on  $B_n$ by some fixed measure on the whole set $\Delta$ or a measure
  not related to any measure on $\Delta$).
  The stratification (\ref{eq:stratification}) and the ensemble of measures
 \begin{equation} \label{eq:measures}
 \{\mu_n\} =  \{\mu_n \mid n \in \mathbb{N}\}
 \end{equation}
   allow one to estimate the asymptotic behavior of
  subsets of $\Delta$. For a subset $R \subseteq \Delta$ the {\em asymptotic
  density} $\rho_\mu(R)$ is defined by the following limit (if it
  exists)
   $$
\rho_\mu(R) = \lim_{n \rightarrow \infty} \mu_n(R \cap B_n).
$$
If  $\mu_n$ is the uniform distribution on the finite set $B_n$ then
   $$  \mu_n(R \cup B_n) = \frac{|R \cap B_n|}{|B_n|}$$
 is the $n$-th {\em frequency}, or probability, to hit an element from $R$ in the ball
 $B_n$. In this case we refer to $\rho_\mu(R)$ as to the {\em asymptotic density} of
 $R$ and  denote it  by $\rho(R)$.

  One can also define the asymptotic densities above using $\lim \sup$
rather then $\lim$, in which event $\rho_\mu(R)$ does  always exist.

We say that a subset $R \subseteq \Delta$ is \emph{generic}   if
$\rho_\mu(R)=1$ and \emph{negligible} if $\rho_\mu(R)  = 0$.  It is
worthwhile to mention  that the  asymptotic densities  not only
 allow one to distinguish between "large" (generic) and "small"
 (negligible) sets, but  give a tool
 to differentiate between various large (or small) sets.  For instance,
 we say that $R$
has asymptotic density $\rho_\mu(R)$ with a {\em super-polynomial
convergence rate} if
$$|\rho_\mu(R) - \mu_n(R \cap B_n)|  = o(n^{-k})$$
for any  $k \in \mathbb{N}$. For brevity,   $R$  is called {\em
strongly generic}  if $\rho_\mu(R) = 1$ with a super-polynomial
convergence rate. The set $R$ is {\em strongly negligible} if its
complement $S - R$ is strongly generic.

  Similarly, one can define exponential convergence rates and
  exponentially generic (negligible) sets.

\subsection{Random subgroups and generating tuples}
  \label{subsec:random-sub-gen}

 In this section we  follow the most commonly used in cryptography
 procedure to generate random subgroups of a given group
  (see for example \cite{AAG}). In brief, the following procedure
is often employed:

\medskip
{\bf  Random Generator of subgroups in $G$}:
 \begin{itemize}
 \item pick a random $k \in \mathbb{N}$ between given boundaries $K_0 \leq k \leq K_1$;
 \item pick randomly $k$ words $w_1, \ldots, w_k \in F(X)$ with fixed length
 range $L_0 \leq |w_i| \leq L_1$;
 \item  output a subgroup $\langle w_1, \ldots, w_k \rangle$ of $G$.
 \end{itemize}
Without loss of generality we may fix from the beginning a single
natural number $k$, instead of choosing it from the finite interval
$[K_0,K_1]$ (by the formula of complete probability the general case
 can be reduced to this one).
Fix $k\in \mathbb{N}$, $k \geq 1$, and focus on the set of all
 $k$-generated subgroups of $G$.

The corresponding descriptions $\delta$, the size function, and the
corresponding stratification of the set of all descriptions  can
be formalized as follows.
By a description $\delta(H)$ of a $k$-generated subgroup $H$ of $G$
 we understand here any $k$-tuple $(w_1, \ldots,w_k)$ of words from
$F(X)$ that generates $H$ in $G$. Hence, in this case
  the space of all descriptions   is
the cartesian product  $F(X)^k$ of $k$ copies of $F(X)$:
  $$\Delta = \Delta_k = F(X)^k.$$
 The size $s(w_1, \ldots, w_k)$
 can be  defined as the total length of the generators
 $$s(w_1, \ldots, w_k) = |w_1| + \ldots + |w_k|,$$
 or as the maximal length of the components:
     $$s(w_1, \ldots, w_k) = \max\{|w_1| , \ldots , |w_k|\}.$$
 Our approach works for both definitions, so  we do not specify which one we use here.
 For $n \in
\mathbb{N}$ denote by $B_n$  the ball of radius $n$ in $\Delta$:
  $$B_n = \{(w_1, \ldots,w_k) \in F(X)^k \mid s(w_1, \ldots,w_k) \leq  n\}.$$
 This gives the required stratification
  $$\Delta = \cup_{n = 1}^{\infty} B_n.$$
 For a subset $M$ of $\Delta$ we define the asymptotic density
 $\rho(M)$ relative to the stratification above assuming the uniform
 distribution on the balls $B_n$:
   $$ \rho(M) = \lim_{n \rightarrow \infty} \frac{|B_n \cap
  M|}{|B_n|}.
  $$
 Notice, that there are several obvious  deficiencies in this
 approach:  we consider subgroups with a fixed number of generators,
  every subgroup may have  distinct $k$-generating tuples,
  every generator can be described by several distinct words from $F(X)$, i.e., our
  descriptions are far from being unique. However, as we have mentioned above, this models
  describe the standard  methods to generate subgroups
  in cryptographic protocols. We refer to \cite{MO} for other approaches.

 \subsection{Asymptotic  properties of subgroups}

Let $G$ be a group with a finite set of generators  $X$ and $k$ a
fixed positive natural number. Denote by $\mathcal{P}$ a property of
descriptions of $k$-generated subgroups of $G$. By $\mathcal{P}(G)$
we denote the set of all descriptions from $\Delta = \Delta_k$ that
satisfy $\mathcal{P}$ in $G$.

\begin{definition}
We say that a property $\mathcal{P} \subseteq \Delta$ of
descriptions of $k$-generated subgroups of $G$ is:
 \begin{itemize}
   \item [1)] {\em asymptotically visible} in $G$ if $\rho(\mathcal{P}(G)) >
   0$;
   \item [2)] {\em generic} in $G$ if   $\rho(\mathcal{P}(G)) = 1$;
    \item [3)] {\em strongly generic} in $G$ if $\rho(\mathcal{P}(G)) =
    1$ and the rate of convergence of $\rho_n(\mathcal{P}(G))$ is super-polynomial;
     \item [4)] {\em exponentially generic} in $G$ if $\rho(\mathcal{P}(G)) =
    1$ and the rate of convergence of $\rho_n(\mathcal{P}(G))$ is
    exponential.

\end{itemize}
\end{definition}

Informally, if $\mathcal{P}$ is asymptotically visible for
$k$-generated subgroups of $G$ then there is a certain  non-zero
probability that a randomly and  uniformly chosen description
$\delta \in \Delta$ of a sufficiently big size results in a subgroup
of $G$ satisfying $\mathcal{P}$. Similarly, if $\mathcal{P}$ is
exponentially generic for $k$-generated subgroups of $G$ then  a
randomly and uniformly chosen description $\delta \in \Delta$ of a
sufficiently big size results in a subgroup of $G$ satisfying
$\mathcal{P}$ with overwhelming probability. Likewise, one can
interpret generic and strongly generic properties of subgroups. If a
set of descriptions $\Delta$ of subgroups of $G$ is fixed, then  we
sometimes abuse the terminology and refer to asymptotic properties
of descriptions of subgroups as asymptotic properties of the
subgroups itself.

\begin{example}
Let $H$ be a  fixed $k$-generated group.  Consider the following
property $\mathcal{P}_H$: a given description $(w_1, \ldots,w_k) \in
F(X)^k$ satisfies $\mathcal{P}_H$ if the subgroup $\langle w_1,
\ldots,w_k \rangle$, generated in $G$ by this tuple, is isomorphic
to $H$. If $\mathcal{P}_H(G)$ is asymptotically visible (generic) in
$\Delta$ then we say that the group $H$ is asymptotically visible
(generic) in $G$ (among $k$-generated subgroups).
 \end{example}

 By  {\em $k$-spectrum}   $Spec_k(G)$ of $G$ we denote the set
of all (up to isomorphism) $k$-generated groups which are
asymptotically visible in $G$.

There are several natural questions about asymptotically visible
subgroups of $G$ that play an important part in cryptography.  For
example, when choosing $k$-generated subgroups of $G$ randomly it
might be useful to know what kind of subgroups you can get with
non-zero probability. Hence the  following question is of interest:
  \begin{problem}
 What is the spectrum $Spec_k(G)$ for a given group $G$ and a natural number $k \geq 1$?
\end{problem}
 More technical,
but also important in applications is the following  question.
  \begin{problem}
  Does the spectrum $Spec_k(G)$ depend on a given finite set of generators
of $G$? \end{problem}
  We will see in due course that answers to these questions play an
  important part in the choice of strong keys
  in some group-based cryptosystems.

\subsection{Groups with generic free basis  property}
\label{subsec:free-basis-property}

 \begin{definition}
 We say that a tuple $(u_1, \ldots, u_k) \in F(X)^k$ has a {\em
free basis} property ($\mathcal{FB}$) in $G$ if it freely generates
a free subgroup in $G$.
  \end{definition}

  In \cite{Toshiaki}  Jitsukawa showed
 that $\mathcal{FB}$  is generic for $k$-generated subgroups of a  finitely generated
 non-abelian  group $F(X)$  for every $k \geq 1$  with respect to the standard basis $X$. Martino, Turner
and  Ventura   strengthened this result in \cite{enric}, they proved
 that $\mathcal{FB}$ is {\em exponentially generic} in  $F(X)$  for
every $k \geq 1$ with respect to the standard basis $X$.
 Recently, it has been
shown in \cite{MO}
 that  $\mathcal{FB}$  is exponentially generic in arbitrary hyperbolic non-elementary
 (in particular, free non-abelian) group for every $k \geq 1$ and  with respect
 to any finite set of generators.

We say that the group $G$ has the {\em generic free basis} property
if $\mathcal{FB}$  is generic in $G$ for every $k \geq 1$ and every
finite generating set of $G$. Similarly, we define groups with {\em
strongly} and {\em exponentially} generic free basis property.
 By $\mathcal{FB}_{gen}$,  $\mathcal{FB}_{st}$, $\mathcal{FB}_{exp}$
  we denote  classes of finitely generated groups with, correspondingly, generic,
  strongly generic, and exponentially generic, free basis property.

The following result gives a host of examples of groups with generic
$\mathcal{FB}$.

\begin{theorem}
\label{th:quotient}
 Let $G$ be a finitely generated group and $N$ a normal subgroup of $G$.
 If the quotient group
$G/N$ is in $\mathcal{FB}_{gen}$, or in $\mathcal{FB}_{st}$, or in
$\mathcal{FB}_{exp}$,  then the whole group $G$ is in the same
class.
\end{theorem}
\begin{proof} Let $H = G/N$ and $\phi:G \rightarrow H$ be the canonical
epimorphism. Fix a finite generating set $X$ of $G$ and a natural
number $k \geq 1$. Clearly, $X^\phi$ is a finite generating set of
$H$. By our assumption, the free basis property is generic in
 $H$  with respect to
the generating set $X^\phi$ and given $k$.  Identifying $x \in X$
with $x^\phi \in H$ we may assume that a  finitely generated
subgroup $A$ of $G$ and the subgroup $A^\phi$ have the same set of
descriptions. Observe now, that for a subgroup $A$ of $G$ generated
by a $k$-tuple $(u_1, \ldots,u_k) \in F(X)^k$ if $A^\phi$ is free
with basis $(u_1^\phi, \ldots,u_k^\phi)$ then $A$ is also free with
basis $(u_1, \ldots,u_k)$. Therefore for each $t \in \mathbb{N}$
$$\frac{|B_t \cap
  \mathcal{FB}(G)|}{|B_t|}   \geq \frac{|B_t \cap
  \mathcal{FB}(H)|}{|B_t|}.$$
 This implies, that if  $\mathcal{FB}(H)$ is generic in $H = G/N$, that $\mathcal{FB}(G)$
 is also generic in $G$, and its  convergence rate  in $G$ is not
 less then the corresponding  convergence rate in $H$,
 as claimed.

\end{proof}

The result above bears on some infinite groups used recently in
group-based cryptography. Braid groups $B_n$ appear as the main
platform in the braid-group cryptography (see
\cite{AAG,Ko1,D4,AAG2}). Recall that the braid group $B_n$ can be
defined by the classical Artin presentation:
\begin{displaymath}
B_n =
 \left\langle
\begin{array}{lcl}\sigma_1,\ldots,\sigma_{n-1} & \bigg{|} &
\begin{array}{ll}
 \sigma_i \sigma_j \sigma_i= \sigma_j \sigma_i \sigma_j & \textrm{if }|i-j|=1 \\
 \sigma_i \sigma_j=\sigma_j \sigma_i & \textrm{if }|i-j|>1
\end{array}
\end{array}
\right\rangle.
\end{displaymath}
Denote by $\sigma_{i,i+1}$ the transposition $(i,i+1)$ in the
symmetric group $\Sigma_n$. The map $\sigma_i \rightarrow
\sigma_{i,i+1}$, $i = 1, \ldots, n$ gives rise to the canonical
epimorphism $\pi:B_n \rightarrow \Sigma_n$. The kernel of $\pi$ is
a subgroup of index $n!$ in $B_n$, termed the {\em pure braid}
group $PB_n$.

\begin{corollary}
 \label{co:pure-braids}
 The free basis property is exponentially generic in the pure braid groups $PB_n$ for $n \geq
 3$.
\end{corollary}
  \begin{proof} It is known (see \cite{Birman}, for example) that
  a pure braid group $PB_n, n\geq 3,$ has the group $PB_3$ as its epimorphic
  quotient, and the group $PB_3$ is isomorphic to $F_2 \times \mathbb{Z}$,
  so $PB_n, n \geq 3,$ has the free group $F_2$ as its quotient.
Now, the result follows from Theorem \ref{th:quotient} and the
strong version of the  Jitsukawa's result \cite{enric,MO,Toshiaki}.
  \end{proof}

 As we have seen a pure braid group $PB_n, n \geq
 3,$ has exponentially generic free basis property and
 it is a subgroup of finite
 index in the braid group $B_n$. However,   at the moment,  we do not have a proof
 that $B_n$ has exponentially generic free basis property. Though,
 we conjecture that this should be true.

\begin{problem}
Is it true that the braid groups $B_n$, $n \geq
 3$, has exponentially generic free basis property?
\end{problem}

In \cite{WM}  partially commutative groups were  proposed as
possible platforms for some cryptosystems. We refer to \cite{Bir}
for more recent discussion on this. By definition a partially
commutative group $G(\Gamma)$ (also called, sometimes, as  right
angled Artin groups, or graph groups, or  trace groups) is a group
associated with a finite graph $\Gamma = (V,E)$, with a set of
vertices $V = \{v_1, \ldots,v_n\}$ and a set of edges $E \subseteq V
\times V$,  by the following presentation:
$$
G(\Gamma) = \langle v_1, \ldots,v_n \mid v_iv_j = v_jv_i \ for  \
(v_i,v_j) \in E \rangle.
$$
Observe, that the group $G(\Gamma)$ is abelian if and only if the
graph $\Gamma$ is complete.
\begin{corollary}
 \label{co:part-comm}
 The free basis property is exponentially generic in
 non-abelian partially commutative groups.

\end{corollary}
 \begin{proof}
   Let $G = G(\Gamma)$ be  a non-abelian partially commutative group
  corresponding to a finite graph $\Gamma$. Then
there are three vertices in $\Gamma$, say $v_1, v_2, v_3$ such that
the complete subgraph $\Gamma_0$ of $\Gamma$ generated by these
vertices is not a triangle. In particular, a partially commutative
group $G_0 = G(\Gamma_0)$ is either a free group $F_3$ (no edges in
$\Gamma_0$), or $(\mathbb{Z} \times \mathbb{Z})  \ast \mathbb{Z}$
(only one edge in $\Gamma_0$), or $F_2 \times \mathbb{Z}$ (precisely
two edges in $\Gamma_0$). Notice that in all three cases the group
$G(\Gamma_0)$ has $F_2$ as its epimorphic quotient. Now, it suffices
to show that $G(\Gamma_0)$ is an epimorphic quotient of $G(\Gamma)$,
which is obtained from $G(\Gamma)$ by adding  to the standard
presentation of $G(\Gamma)$ all the relations of the type $v = 1$,
where $v$ is a vertex of $\Gamma$ different from $v_1, v_2, v_3$.
This shows that $F_2$ is a quotient of $G(\Gamma)$ and the result
follows from Theorem \ref{th:quotient}.
  \end{proof}

Observe, that some other groups, that have been proposed as
platforms in based-group cryptography, do not have non-abelian free
subgroups at all, so they do not have free basis property for $k
\geq 2$. For instance,  in \cite{P} the Grigorchuk groups were used
as a platform. Since these   groups  are periodic (i.e., every
element has finite order) they  do not contain non-trivial free
subgroups. It is not clear what are asymptotically visible subgroups
in Grigorchuk groups.  As another example, notice that in \cite{SU}
authors put forth the Thompson group $F$ as a platform. It is known
that there are no non-abelian free subgroups in $F$ (see, for
example, \cite{BS}), so $F$ does not have free basis property.
Recently, some interesting results were obtain on the spectrum
$Spec_k(F)$  in \cite{Elder}.

\subsection{Quasi-isometrically embedded subgroups}

In this section we discuss another property of subgroups of $G$ that
plays an important part in our cryptanalysis of group based
cryptosystems.

Let $G$ be a group with a finite generating set $X$. The {\em Cayley
graph} $\Gamma(G,X)$ is an $X$-labeled directed  graph with the
vertex set $G$ and such that any two vertices $g,h \in G$ are
connected by an edge from $g$ to $h$ with a label $x \in X$ if and
only if  $gx = h$ in $G$. For convenience we usually assume that the
set $X$ is closed under inversion, i.e., $x^{-1} \in X$ for every $x
\in X$. One can introduce a metric $d_X$ on $G$ setting $d_X(g,h)$
equal to the length of a shortest word in $X^{\pm 1} = X \cup
X^{-1}$ representing the element $g^{-1}h$ in $G$. It is easy to see
that $d_X(g,h)$ is equal to the length of a shortest path from $g$
to $h$ in the Cayley graph $\Gamma(G,X)$. This turns $G$ into a
metric space $(G,d_X)$. By $l_X(g)$ we denote the length of a
shortest word in generators $X^{\pm 1}$ representing the element
$g$, clearly $l_X(g) = d_X(1,g)$.

  Let $H$ be a subgroup of $G$ generated by a finite set of elements
  $Y$. Then there are two metrics on $H$:  the first
  one is $d_Y$ described above and the other one is the metric
  $d_X$ induced from the metric space $(G,d_X)$ on the subspace
  $H$. The following notion allows one to compare these metrics.
    Recall that a map $f:M_1 \rightarrow M_2$ between two metric spaces
    $(M_1,d_1)$ and $(M_2,d_2)$ is a {\em quasi-isometric embedding} if there are
    constants $\lambda > 1, c > 0$ such that for every elements $x,y
    \in M_1$ the following inequalities hold:
      \begin{equation} \label{eq:quasi-isom}
       \frac{1}{\lambda} d_1(x,y) - c \leq  d_2(f(x),f(y)) \leq
      \lambda d_1(x,y) + c.
       \end{equation}
 In particular, we say that a subgroup $H$ with a finite set of
 generators $Y$ is {\em quasi-isometrically embedded} into $G$ if
 the inclusion map $i:H \hookrightarrow G$ is a quasi-isometric
 embedding $i:(H,d_Y) \rightarrow
 (G,d_X)$.  Notice, that in this case the right-hand inequality in
 (\ref{eq:quasi-isom}) always holds, since for all $f, h \in H$
  $$d_X(i(f),i(h)) \leq       \max_{y \in Y}\{l_X(y)\} \cdot d_Y(f,h).$$
  Therefore, the definition of quasi-isometrically embedded subgroup
  takes the following simple form (in the notation above).
    \begin{definition}
    Let $G$ be a group with a finite generating set $X$ and
  $H$  a subgroup of $G$ generated by a finite set of elements
  $Y$. Then $H$ is
  {\em quasi-isometrically embedded} into $G$ if  there are
    constants $\lambda > 1, c > 0$ such that for every elements $f,h
    \in H$ the following inequality holds:
  \begin{equation} \label{eq:quasi-isom-short}
       \frac{1}{\lambda} d_Y(f,h) - c \leq  d_X(f,h).
       \end{equation}
\end{definition}

   It follows immediately from the definition, that if $X$ and $X^\prime$ are two finite generating sets of
 $G$ then the metric spaces $(G,d_X)$ and $(G,d_{X^\prime})$ are quasi-isometrically
 embedded into each other. This implies  that the notion of quasi-isometrically embedded subgroups
 is independent of the choice of finite generating sets in $H$ or in $G$
 (though the constants $\lambda$ and $c$ could be different).

\begin{definition}
Let $G$ be a group with a finite generating set $X$. We say that a
tuple $(u_1, \ldots, u_k) \in F(X)^k$ has a $\mathcal{QI}$
(quasi-isometric embedding) property in $G$ if the subgroup it
generates in $G$ is quasi-isometrically embedded into $G$.
\end{definition}
Denote by $\mathcal{QI}(G)$ the set of all tuples in $F(X)^k$ that
satisfy the $\mathcal{QI}$ property in $G$. We term  the property
$\mathcal{QI}$ is {\em generic} in $G$ if $\mathcal{QI}(G)$ is
generic in $G$  for every $k \geq 1$ and every finite generating set
of $G$. Similarly, we define groups with {\em strongly} and {\em
exponentially} generic quasi-isometric embedding subgroup property.
Denote by $\mathcal{QI}_{gen}$, $\mathcal{QI}_{st}$,
$\mathcal{QI}_{exp}$
  classes of finitely generated groups with, correspondingly, generic,
  strongly generic, and exponentially generic, quasi-isometric embedding
  subgroup property.

  It is not hard to see that {\em every}  finitely generated subgroup of
  a finitely generated
  free group $F$ is quasi-isometrically embedded in $F$, so $F \in \mathcal{QI}_{exp}$.

  The following
  result gives further examples of groups with quasi-isometric embedding
  subgroup property.

   Let $G \in \mathcal{FB}_{gen} \cap \mathcal{QI}_{gen}$.
Notice,  that the intersection of two generic
  sets $\mathcal{FB}(G) \subseteq F(X)^k$ and $\mathcal{QI}(G) \subseteq F(X)^k$ is again
  a generic set  in $F(X)^k$,  so the set $\mathcal{FB}(G) \cap \mathcal{QI}(G)$  of all
  descriptions $(u_1, \ldots,u_k) \in F(X)^k$  that freely generate a quasi-isometrically embedded subgroup
  of $G$,  is generic in $F(X)^k$.   Observe, that by the remark above, and the result on free basis property in free groups,
    $\mathcal{FB}_{gen} \cap \mathcal{QI}_{gen}$  contains all free
  groups of finite rank. The argument applies also to the
  strongly generic and exponentially generic variations of the
  properties. To unify references  we will use the following
  notation: $\mathcal{FB}_{\ast} \cap \mathcal{QI}_{\ast}$ for $\ast
  \in \{gen, st, exp\}$.

    \begin{theorem}
    \label{th:NQ}
       Let $G$ be a finitely generated group with  a
       quotient $G/N$. If $G/N \in \mathcal{FB}_{\ast} \cap \mathcal{QI}_{\ast}$
       then $G \in  \mathcal{FB}_{\ast} \cap \mathcal{QI}_{\ast}$  for any $\ast
  \in \{gen, st, exp\}$.
    \end{theorem}
  \begin{proof}
  Let $G$ be a finitely generated group generated by $X$, $N$ a normal subgroup of
  $G$ such that the quotient
  $G/N$ is in $\mathcal{FB}_{\ast} \cap \mathcal{QI}_{\ast}$. Let  $\phi:G \rightarrow G/N$ be
  the canonical epimorphism. By Theorem \ref{th:quotient}
  $G \in  \mathcal{FB}_{\ast}$, so it suffices to  show now that $G \in \mathcal{QI}_{\ast}$.

   Let $H$ be a $k$-generated
  subgroup with a  set of generators $Y = (u_1, \ldots,u_k) \in F(X)^k$.
  Suppose that $Y \in \mathcal{FB}_{\ast}(G/N) \cap \mathcal{QI}_{\ast}(G/N)$,
   i.e., the image $Y^\phi$ of $Y$ in $G/N$
  freely generates a free group quasi-isometrically embedded into $G/N$.
   Observe, first, that for every element
  $g \in G$  one has $l_X(g) \geq l_{X^\phi}(g^\phi)$, where
  $l_{X^\phi}$  is the length  on $G/N$ relative to the set
  of generators $X^\phi$. Since the subgroup $H^\phi$ is quasi-isometrically embedded
  into $G/N$ the metric space $(H^\phi,d_{Y^\phi})$ quasi-isometrically  embeds into
   $(G^\phi,d_{X^\phi})$. On the other hand, $\phi$ maps  the subgroup $H$ onto the subgroup
   $H^\phi$ isomorphically (since both are free groups with the corresponding bases),
   such that for any $h \in H$ $d_Y(h) =
   d_{Y^\phi}(h^\phi)$. Now we can deduce the following inequalities
   for $g,h \in H$:
    $$   \frac{1}{\lambda}d_{Y}(g,h) - c = \frac{1}{\lambda}d_{Y^\phi}(g^\phi,h^\phi) - c
     \leq d_{X^\phi}(g^\phi,h^\phi) \leq
    d_X(g,h)$$
 where $\lambda$ and $c$ come from the quasi-isometric embedding of
 $H^\phi$ into $G/N$. This shows that  $H$ is quasi-isometrically embedded into $G$,  as required.

     \end{proof}

\begin{corollary}
 \label{co:NQ-pure-braids}
 The following groups are in $\mathcal{FB}_{exp} \cap \mathcal{QI}_{exp}$:
  \begin{enumerate}
   \item [1)] Pure braid groups $PB_n$, $n \geq 3$;
  \item [2)] Non-abelian partially commutative groups $G(\Gamma)$.
   \end{enumerate}
\end{corollary}
\begin{proof} The arguments in Corollaries \ref{co:pure-braids}, \ref{co:part-comm}
show that the groups $PB_n$, $n \geq 3$,  and  $G(\Gamma)$,
non-commutative, have quotient isomorphic to  the free group $F_2$.
Now the result follows from Theorems \ref{th:quotient} and
\ref{th:NQ}.
\end{proof}

\section{Anshel-Anshel-Goldfeld scheme}
\label{se:as-dominance-crypto}

In this section we discuss  the Anshel-Anshel-Goldfeld (AAG)
cryptosystem for public key exchange \cite{AAG} and touch briefly on
its  algorithmic security.

\subsection{Description of Anshel-Anshel-Goldfeld scheme}
 \label{subsec:AAG-scheme}

Here we give a general description of  the Anshel-Anshel-Goldfeld
cryptosystem.

Let $G$ be a  group with a finite generating set $X$, it is called
the {\em platform} of the scheme. We assume that elements $w$ in $G$
have unique normal forms $\bar w$ such that it is "hard" to
reconstruct $w$ from $\bar w$ and there is a "fast" algorithm to
compute $\bar w$ when given $w$.  We do not discuss here the
security issues of these two components of the platform $G$, leaving
this for the future.

The Anshel-Anshel-Goldfeld key exchange protocol requires the
following sequence of steps. Alice [Bob resp.] chooses a random
subgroup  of $G$
    $$A = \langle  a_1, \ldots, a_m \rangle \ \ \ [B = \langle b_1, \ldots, b_n \rangle \mbox{ resp.}]$$
by randomly choosing generators $a_1, \ldots, a_m$ [$b_1, \ldots,
b_n$ resp.] as words in $X^{\pm 1}$, and makes it public. Then Alice
[Bob resp.] chooses randomly a secret element $a = u(a_1, \ldots,
a_m) \in A$ [$b = v(b_1, \ldots, b_n) \in B$ resp.] as a product of
the generators of $A$ [$B$ resp.] and their inverses, takes the
conjugates $b_1^a, \ldots, b_n^a$ [$a_1^b, \ldots, a_m^b$ resp.],
encodes them by taking their normal forms $\overline{b_i^a}$
[$\overline{a_j^b}$ resp.], and makes these normal forms public:
    $$\overline{b_1^a}, \ldots, \overline{b_n^a} \ \ \ [\overline{a_1^b}, \ldots, \overline{a_m^b} \mbox{ resp.}].$$
Afterward, they both can compute the secret shared element of $G$:
 $$a^{-1}a^b = [a,b] = (b^a)^{-1}b$$
and take its normal form as the secret shared key.

\subsection{Security assumptions of AAG scheme}
 \label{subsec:sec-ass-AAG}

In this section we briefly discuss computational security features
of the AAG cryptosystem. Unfortunately, in the original description
of AAG the authors did not state  precisely what are the security
assumptions  that make the  system difficult to break. Here we dwell
on several possible assumptions of this type, that  often occur,
though sometimes implicitly,  in the literature on the AAG
cryptosystem.

 It seems that the security of AAG relies on
the computational hardness of the following, relatively new,
computational problem in group theory:

\bigskip
\noindent
 {\bf  AAG Problem:}  given the whole public information
from the scheme AAG, i.e., the group $G$, the elements $ a_1,
\ldots, a_m $, $b_1, \ldots, b_n$, and  $\overline{b_1^a}, \ldots,
\overline{b_n^a}, \overline{a_1^b}, \ldots, \overline{a_n^b}$ in a
group $G$, find the shared secret key $[a,b]$.

\bigskip
This problem is not a standard   group-theoretic problem, not much
is known about its complexity, and it is quite technical to
formulate. So it would be convenient  to reduce this problem to a
standard algorithmic problem in groups or to a combination of such
problems. The following problems seem to be relevant here and they
attracted quite a lot of attention recently, especially in the braid
groups -- the original platform for AAG \cite{AAG}. We refer to
papers \cite{BH}, \cite{BGM1}, \cite{BGM2}, \cite{G}, \cite{L},
\cite{LL}. Nevertheless, the precise relationship between these
problems and AAG  is unclear, see \cite{SU2} for more details.

\bigskip
\noindent {\bf The  Conjugacy Search  Problem (CSP):}   given $u, v
\in G$ such that an equation $u^x = v$ has a solution in $G$, find a
solution.

\bigskip
\noindent {\bf  The Simultaneous  Conjugacy Search  Problem (SCSP):}
given $u_i, v_i \in G$,  such that a system  $u_i^x = v_i$, $i = 1,
\ldots, n$ has a solution in $G$,  find a solution.

\bigskip
\noindent
 {\bf The Simultaneous  Conjugacy Search  Problem relative to a subgroup (SCSP*):} given
$u_i, v_i \in G$  and a finitely generated subgroup $A$ of $G$ such
that  a system $u_i^x = v_i$, $i = 1, \ldots, n$ has a solution in
$A$,  find such a solution.

\begin{remark}
Observe, that if the Word Problem is decidable in $G$ then all
the problems above are also decidable. Indeed, one can enumerate all
possible elements $x$ (either in $G$ or in the subgroup $A$),
substitute them one-by -one into the equations, and check, using the
decision algorithm for the Word Problem in $G$,  if $x$ is a
solution or not. Since the systems above have some  solutions
 this algorithm will eventually find one. However, the main problem
here is not about decidability, the problem  is whether or not one
can find a solution sufficiently "quickly", say in polynomial time
in the size of the inputs.
 \end{remark}

The following result is easy.

\begin{lemma}
\label{le:AAG-reduced}
 For any group $G$ the AAG problem can be
reduced in linear time to the problem SCSP*.
\end{lemma}
 \begin{proof} Suppose in a finitely generated group $G$
  we are given the public data from the AAG scheme, i.e.,
  the   subgroups
  $$A = \langle  a_1, \ldots, a_m \rangle,  \ \ B = \langle b_1,
\ldots, b_n \rangle, $$
 and the elements
  $\bar{b_1^a}, \ldots, \bar{b_n^a}$ and $\bar{a_1^b},
\ldots, \bar{a_n^b}$. If the problem SCSP relative to subgroups $A$
and $B$ is decidable in $G$,  then solving a system of equations
  \begin{equation}
  \label{eq:b}
   b_1^x=\bar{b_1^a}, \ldots , b_n^x =\bar{b_n^a}
 \end{equation}
in $A$  one can find a solution $u \in A$.
 Similarly, solving a system of equations
 \begin{equation}
  \label{eq:a} a_1^y =\bar{a_1^b}, \ldots , a_m^y =\bar{a_m^b}
   \end{equation}
 in $B$ one can find a solution $v \in B$. Notice, that all  solutions
 of the system (\ref{eq:b})  are elements of the form  $ca$ where $c$
 is an arbitrary element from the centralizer $C_G(B)$, and all
  solutions of the system (\ref{eq:a}) are of the form $db$ for some  $d
  \in C_G(A)$. In this case, obviously $[u,v] = [ca,db]= [a,b]$ gives a solution to the
   the AAG problem.
\end{proof}

Clearly, in some groups, for example, in abelian  groups  AAG
problem as well as the SCSP* are both decidable in polynomial time,
which makes them (formally) polynomial time equivalent. We will see
in Section \ref{subsec:LBA-free-groups} that   SCSP* is easy in free
groups.

It is not clear, in general,  whether the SCSP is any harder or
easier than the CSP. In hyperbolic groups SCSP, as well as CSP,  is
easy \cite{BH}.

There are indications that in  finite simple  groups, at least
generically, the SCSP* is not harder than the standard CSP (since,
in this case,  two randomly chosen elements generate the whole
group). We refer to a preprint \cite{GMMU} for a brief discussion on
complexity of these problems.

It is interesting to get some information on the following problems,
which would shed some light on the complexity of AAG problem.

\begin{problem}
\begin{itemize}
 \item [1)] In which groups AAG problem is poly-time equivalent to the
SCSP*?
 \item [2)] In which groups SCSP*   is harder than the SCSP?
 \item [3)] In which groups SCSP is harder (easier) than CSP?
   \end{itemize}
\end{problem}

In the rest of the paper we study the hardness of SCSP* in various
groups and analyze some of the most successful attacks on AAG from
the view-point of asymptotic mathematics.

\section{ Length Based Attacks}
 \label{subse:LBA}

The intuitive idea of the length based attack (LBA) was first put in
the paper \cite{HT} by J. Hughes and A. Tannenbaum. Later it was
further developed in a joint paper \cite{GKTTV} by Garber, Kaplan,
Teicher, Tsaban, and Vishne where the authors gave an experimental
results concerning the success probability of LBA that suggested
that very large computational power is required for this method to
successfully solve the Conjugacy Search Problem.

Recently, the most successful variation of this attack for braid
groups was developed in \cite{MU3} where the authors suggested to
use a heuristic algorithm for approximation of the geodesic length
of braids in conjunction with LBA. Furthermore, the authors analyzed
the reasons for success/failure of their variation of the attack, in
particular the practical importance of Alice's and Bob's subgroups
$A$ and $B$ being non isometrically embedded and being able to
choose the elements of these subgroups distorted in the group (they
refer to such elements as peaks).

In this section we rigorously prove that the same results can be
observed in much  larger classes of groups. In particular our
analysis works for the class $\mathcal{FB}_{exp}$ and, hence, for
free groups, pure braid groups, locally commutative non-abelian
groups, etc.

\subsection{A general description}
  \label{subsec:general-LBA}

Since LBA is an attack on AAG scheme the inputs for LBA are
precisely the inputs for AAG algorithmic problem. Moreover, in all
its variations LBA attacks AAG via solving the corresponding
conjugacy equations given in a particular instance of AAG. In what
follows we take a slightly more general approach and view the length
based attack (LBA) as a correct partial search deterministic
algorithm of a particular type for the Simultaneous Conjugacy Search
Problem relative to a subgroup in a given group $G$. In this case
LBA is employed to solve SCSP*, not AAG. Below we describe a basic
LBA in its most simplistic form.

Let $G$ be a group with a finite generating  set $X$. Suppose we are
given a particular instance of the SCSP*, i.e., a system  of
conjugacy equations $u_i^x = v_i, i= 1, \ldots, m$ which has a
solution in a subgroup $A = \langle Y  \rangle$ generated by a
finite set $Y$ of elements in  $G$ (given by words in $F(X)$). The
task is to find such a solution in $A$.  The main idea of LBA is
very simple and it is based on the following assumptions:
\begin{enumerate}
    \item [\bf(L1)]
for arbitrary "randomly chosen" elements $u,w \in G$ one has
$l_X(u^w) > l_X(u)$;
    \item[\bf(L2)]
for "randomly chosen" elements $w, y_1, \ldots,  y_k$ in $G$ the
element $w$ has minimal $l_X$-length among all elements of the type
$w^y$, where $y$ runs over the subgroup of $G$ generated by $y_1,
\ldots,  y_k$.
\end{enumerate}
   It is not obvious at all  whether this assumption is
realistic or not, or even  how to formulate it correctly. We will
return to these  issues  in due course. Meantime, to make use of the
assumptions above we assume that we are given an algorithm
$\mathcal{A}$ to compute the length function $l_X(w)$  for a given
element $w \in G$.

Consider Alice' public conjugates $\bar{b}_1^a, \ldots,
\bar{b}_n^a$, where $a = a_{s_1}^{\varepsilon_1} \ldots
a_{s_L}^{\varepsilon_L}$. Essentially each $\bar{b}_i^a$ is a result
of a sequence of conjugations of $b_i$ by the factors of $A$:
\begin{equation}\label{eq:seq_conj}
\begin{array}{rrcl}
& & b_i & \\
& & \downarrow & \\
& a_{s_1}^{-\varepsilon_1} & b_i & a_{s_1}^{\varepsilon_1} \\
& & \downarrow & \\
& a_{s_2}^{-\varepsilon_2} a_{s_1}^{-\varepsilon_1} & b_i & a_{s_1}^{\varepsilon_1} a_{s_2}^{\varepsilon_2} \\
& & \downarrow & \\
& & \ldots & \\
& & \downarrow & \\
\bar{b}_i^a = & a_{s_{L}}^{-\varepsilon_{L}} \ldots a_{s_2}^{-\varepsilon_2}
a_{s_1}^{-\varepsilon_1} & b_i & a_{s_1}^{\varepsilon_1} a_{s_2}^{\varepsilon_2}
 \ldots a_{s_{L}}^{\varepsilon_{L}}\\
\end{array}
\end{equation}
A conjugating sequence is the same for each $b_i$ and is defined by
the private key $a$. The main goal of the attack is to reverse the
sequence (\ref{eq:seq_conj}) and going back from the bottom to the
top recover each conjugating factor.  If successful the procedure
will result in the actual conjugator as a product of elements from
$\la$.

The next algorithm is the simplest realization of LBA called the
best descend LBA. It takes as an input three tuples
$(a_1,\ldots,a_{m})$, $(b_1,\ldots,b_{n})$, and $(c_1, \ldots, c_n)$
where the last tuple is assumed to be $\bar{b}_1^a, \ldots,
\bar{b}_n^a$. The algorithm is a sequence of the following steps:
\begin{itemize}
    \item[$-$]
{\bf (Initialization)} Put $x = \varepsilon$.
    \item[$-$]
{\bf (Main loop)} For each $i=1, \ldots, n$ and $\varepsilon = \pm
1$ compute $l_{i, \varepsilon} = \sum_{j=1}^n l_X(a_i^{-\varepsilon}
c_j a_i^{\varepsilon})$.
\begin{itemize}
    \item
If for each $i=1, \ldots, n$ and $\varepsilon = \pm 1$ the
inequality $l_{i, \varepsilon} > \sum_{j=1}^n l_X(c_j)$ is satisfied
then output $x$.
    \item
Otherwise pick $i$ and $\varepsilon$ giving a least value
$l_{i,\varepsilon}$. Multiply $x$ on the right by
$a_i^{\varepsilon}$. For each $j = 1, \ldots, n$ conjugate $c_j =
a_i^{-\varepsilon} c_j a_i^{\varepsilon}$. Continue.
\end{itemize}
    \item[$-$]
{\bf (Last step)} If $c_j = b_j$ for each $j = 1, \ldots, n$ then
output the obtained element $x$. Otherwise output $Failure$.
\end{itemize}
Other variations of LBA suggested in \cite{MU3} are LBA with
Backtracking and Generalized LBA. We refer to \cite{MU3} for a
detailed discussion on this.

One can notice that instead of the length function $l_X$ one can use
any other objective function satisfying assumptions (L1) and (L2).
In this work besides $l_X$ we analyze the behavior of modifications
of LBA relative to the following functions:
\begin{itemize}

    \item[\bf(M1)]
Instead of computing the geodesic length $l_X(v_i)$ of the element
$v_i \in G$ compute the geodesic length $l_Z(v_i)$ in the subgroup
$H$ generated by $Z = \{u\} \cup Y$ (clearly, $v_i \in H$). In this
case, LBA in $G$ is reduced to LBA in $H$, which might be easier.
 We term $l_Z$ the {\em inner} length in LBA.
    \item[\bf(M2)]
It might be difficult to compute the lengths $l_X(w)$ or $l_Z(w)$.
In this case, one can try to compute some "good", say linear,
approximations of $l_X(w)$ or $l_Z(w)$, and then use some
heuristic algorithms to carry over LBA (see \cite{MU3}).
\end{itemize}
These modifications can make LBA much more efficient as we will see
in the sequel.

In what follows our main interest is in the generic time complexity
of LBA. To formulate this precisely one needs to describe the set of
inputs for LBA and the corresponding distribution on them.

Recall that an input for SCSP* in a given  group $G$ with a fixed
finite generating set $X$ consists of a finitely generated subgroup
$A = \langle a_1, \ldots, a_k \rangle $ of $G$ given by a $k$-tuple
$(a_1, \ldots, a_k) \in F(X)^k$, and a finite system of conjugacy
equations  $u_i^x = v_i$, where $u_i, v_i \in F(X)$,  $i = 1,
\ldots, m$, that has a solution in $A$. We denote this data by
 $\alpha = (T,b)$, where $T =(a_1, \ldots, a_k, u_1, \ldots, u_m)$
 and $b = (v_1, \ldots,v_m)$. The distinction that we make here between $T$ and
 $b$ will be in use later on. For fixed positive integers $m, k$ we
 denote the set of all inputs $\alpha = (T,b)$ as above by
 $I_{k,m}$.

The standard procedure to generate a "random" input of this type in
AAG protocol is as follows.

\medskip
{\bf  A Random Generator of inputs for LBA in a given $G$}:
 \begin{itemize}
 \item pick a random $k \in \mathbb{N}$ from a fixed interval  $K_0 \leq k \leq K_1$;
 \item pick randomly $k$ words $a_1, \ldots, a_k \in F(X)$ with the length in fixed interval
  $L_0 \leq |w_i| \leq L_1$;
  \item pick a random $m \in \mathbb{N}$ from a fixed interval $M_0
  \leq m \leq M_1$;
  \item pick randomly $m$ words $u_1, \ldots, u_m \in F(X)$ with the length in fixed interval
  $N_0 \leq |u_i| \leq N_1$;
 \item  pick a random element $w$ from the subgroup $A = \langle a_1, \ldots, a_k \rangle $,  as a random
 product $w = a_{i_1}a_{i_2} \ldots a_{i_c}$ of elements from $\{a_1, \ldots, a_k\}$
 with the number of factors  $c$ in a fixed interval $P_1 \leq c \leq P_2$;
 \item conjugate  $v_i = u_i^w$ and compute the normal form $\tilde{v}_i$ of
 $v_i$, $i = 1, \ldots, m$.
 \end{itemize}

\medskip

As we have argued in Section \ref{subsec:random-sub-gen} one can fix
the numbers $k, m$, and the number of factors $c$ in the product
$w$,  in advance. Observe, that the choice of the elements $v_1,
\ldots, v_m$ is completely determined by the choice of the tuple
 $T = (a_1, \ldots, a_k,u_1, \ldots, u_m) \in F(X)^{k+m}$ and the word $w$.

 Notice, that the  distribution on the subgroups
  $H = \langle T\rangle$ (more precisely,
  their descriptions from $F(X)^{k+m}$) that comes from the random
 generator above coincides with the distribution on the $(k+m)$-generated
 subgroups (their descriptions) that was described in Section \ref{subsec:random-sub-gen}.
 We  summarize this in the following remark.

\begin{remark}\label{re:random-LBA}
\begin{itemize}
 \item [1)] The choice of a tuple  $T = (a_1, \ldots, a_k,u_1, \ldots,
u_m) \in F(X)^{k+m}$ precisely corresponds to the choice of
generators of random subgroups described in Section
\ref{subsec:random-sub-gen}.
 \item [2)] Asymptotic properties of the subgroups generated by $T$
  precisely correspond to the  asymptotic properties
 of subgroups discussed in Section \ref{se:asym-properties}.
 \end{itemize}
\end{remark}

\subsection{LBA in free groups}
  \label{subsec:LBA-free-groups}

In this section we discuss LBA in free groups. It is worthwhile to
mention here that there are fast (quadratic time) algorithms  to
solve SCSP* and, hence, AAG in free groups (see Section
\ref{subsec:SCSP-free-groups}). However, results on LBA in free
groups will serve us as a base for solving SCSP* in many other
groups.

Let $k$ be a fixed positive natural number.
 We say that cancelation in a set of words $Y = \{y_1, \ldots,y_k\}
 \subseteq
 F(X)^k$ is at most $\lambda$, where $\lambda \in (0,1/2)$,  if for any $u,v \in Y^{\pm 1}$
  the amount of cancelation in the
 product $uv$ is strictly less then $\lambda \min\{l_X(u),
 l_X(v)\}$, provided $u \neq v^{-1}$ in $F(X)$.

\begin{lemma} \label{le:small_cond}
If the set $Y = \{y_1, \ldots,y_k\}$ satisfies $\lambda$-condition
for some $\lambda \in (0,1/2)$ then:
\begin{itemize}
    \item
The set $Y$ is Nielsen reduced. In particular, $Y$ freely generates a free subgroup
and any element $w \in \gp{Y}$ can be uniquely represented as a reduced word in the generators $Y$
and their inverses.
    \item
The Membership Search Problem for a subgroup $\gp{Y}$ (see  Section
\ref{subsec:MP-free-groups} for details) is decidable in linear
time.
    \item
The geodesic length for elements of a subgroup $\gp{Y}$ (see Section
\ref{subsec:algor-problems-geodesic} for details) is computable in
linear time.
\end{itemize}
\end{lemma}

\begin{proof}
Easy exercise.

\end{proof}

Moreover, the following result is proved in \cite{enric}.

\begin{theorem}\label{th:small_canc}
Let $\lambda \in (0,1/2)$.
The set $S$ of $k$-tuples $(u_1, \ldots, u_k) \in F(X)^{k}$
satisfying $\lambda$-condition is exponentially generic and,
hence, the set of $k$-tuples which are the Nielsen reduced in $F(X)$ is exponentially generic.
\end{theorem}

Now we are ready to discuss the generic complexity of LBA in free
groups.

\begin{theorem}\label{th:LBA-free-groups}
Let $F(X)$ be a free group with basis $X$. Then LBA with respect to
the inner length $l_Z$ solves SCSP* in linear time on an
exponentially generic set of inputs.
\end{theorem}

\begin{proof}
Let $n$ and $m$ be fixed positive integers.
 Denote by $S$ a set
 of $(n+m)$-tuples $(u_1, \ldots, u_n, a_1, \ldots,
a_m) \in F(X)^{n+m}$ that satisfy $1/4$-condition.
It follows from Theorem \ref{th:small_canc} that the set $S$
is exponentially generic.

Furthermore, the system of conjugacy equations associated with such
a tuple $Z = (u_1, \ldots, u_n, a_1, \ldots, a_m)$ has the form
$$
\left\{
\begin{array}{l}
v_1 = u_1^x \\
\ldots \\
v_n = u_n^x, \\
\end{array}
\right.
$$
where $v_i$ belong to the subgroup  $\langle Z \rangle$ generated by
$Z$ and $x$ is searched in the same subgroup. By Lemma
\ref{le:small_cond} one can find expressions for $v_i$ in terms of
the generators $Z$ in linear time.  Now, since the generators $a_1,
\ldots, a_m$ are part of the basis of the subgroup $\langle Z
\rangle$ it follows that LBA relative to $l_Z$ successfully finds a
solution $x= w(a_1, \ldots, a_m)$ in linear time.

 \end{proof}

\subsection{LBA in groups from $\mathcal{FB}_{exp}$}
 \label{subsec:LBA-FB}

The result above for free groups is not very surprising because of
the nature of cancelation in free groups. What, indeed, looks
surprising is that  LBA works generically
 in some other groups  which seem to be very different from free groups.
 In this and the next section we outline a general mathematical explanation why  LBA
has a high rate of success in various groups, including the braid
groups. In particular, it will be clear why Modification (M1) of
LBA, which was discussed in Section \ref{subsec:general-LBA}, is
very robust, provided  one can compute the geodesic length in
subgroups.

We start with a slight generalization of the result of Theorem
\ref{th:LBA-free-groups}. Recall (from Section
\ref{subsec:general-LBA}) that inputs for  LBA, as well as for
SCSP*, can be described in the form $\alpha = (T,b)$,  where $T =
(a_1, \ldots, a_k,u_1, \ldots, u_m) \in F(X)^{k+m}$ and $b = (v_1,
\ldots, v_m)$, such that there is a solution of the system $u_i^x =
v_i$ in the subgroup $A = \langle a_1, \ldots, a_k \rangle$.

\begin{lemma}
\label{le:I-free}
 Let $G$ be a  group with a finite generating set
$X$ and $I_{k,m}$ a set of all inputs $(T,b)$ for LBA in $G$. Put
 $$I_{free} = \{(T,b) \in I_{k,m} \mid T \ freely \ generates \ a \ free \ subgroup \ in \ G\}.$$
 Suppose there is an exponentially generic subset $S$ of $I_{free}$ and an algorithm  $\mathcal{A}$ that
  computes the geodesic length $l_T$ of elements from the subgroup $\langle T
  \rangle$, $(T,b) \in  S$,
 when these elements are given as words from $F(X)$.
 Then there is an exponentially generic subset $S'$ of $I_{free}$ such that on inputs
 from $S'$ LBA halts and outputs a solution for the related SCSP* in at most
 quadratic time relative to  the algorithm $\mathcal{A}$.
\end{lemma}
  \begin{proof} The result directly follows from Theorem
  \ref{th:LBA-free-groups}.
 \end{proof}

 Let $G \in \mathcal{FB}_{exp}$. In the next theorem we prove that
 the time complexity of SCSP*
 on an exponentially generic set of inputs is at most quadratic relative
 to the time complexity of the problem of computing the geodesic length in finitely
 generated subgroup of $G$.

\begin{theorem}
\label{th:LBA-FB}
 {\bf (Reducibility to subgroup-length function)}
 Let $G$ be a  group with exponentially generic free
basis property and $X$ a finite generating set of $G$. Then there is
an exponentially generic subset $S$ of the set $I_{k,m}$  of all
inputs for LBA in $G$ such that on inputs
 from $S$ LBA relative to $l_T$ halts and outputs a solution for the related SCSP*. Moreover,
 the time complexity   of LBA on inputs from
 $S$ is at most quadratic relative to  the algorithm $\mathcal{A}$ that
  computes the geodesic length $l_T$ of elements from the subgroup $\langle T \rangle$
 when these elements are given as words from $F(X)$.
\end{theorem}
\begin{proof}
   By Lemma \ref{le:I-free} there is an exponentially  generic  subset $S$ of
   $I_{free}$ such that on inputs
 from $S$ LBA halts and outputs a solution for the related SCSP*. Moreover,
 the time complexity   of LBA on inputs from
 $S$ is at most quadratic relative to  the algorithm $\mathcal{A}$ that
  computes the geodesic length $l_T$ of elements from the subgroup $\langle T \rangle$
 when these elements are given as words from $F(X)$. It suffices to
 show now that the set $I_{free}$ is exponentially generic in the
 set of all inputs $I$ for LBA in $G$. By Remark \ref{re:random-LBA}
 asymptotic density of the set $I_{free}$ in $I$ is the same as the
 asymptotic density of the set of tuples $T \in F(X)^{k+m}$ which
 have free basis property in $G$. Since $G$ is in
 $\mathcal{FB}_{exp}$ this set is exponentially generic in
 $F(X)^{k+m}$, so is $I_{free}$ in $I$. This proves the theorem.
\end{proof}

\section{Computing the geodesic length in a subgroup}
 \label{subsec:geodesic-length-general}

For groups $G \in \mathcal{FB}_{exp}$ Theorem \ref{th:LBA-FB}
reduces in quadratic time  the time complexity of LBA on an
exponentially generic set of inputs to the time complexity  of the
problem of computing the geodesic length in finitely
 generated subgroups of $G$. In this section we discuss time complexity
 of  algorithms  to compute the geodesic length
 in a  subgroup of $G$. This discussion is related to Modification 2 of LBA,
 introduced in Section \ref{subsec:general-LBA}.  In particular, we
 focus on the situation when we do not have fast algorithms to
 compute the geodesic length of elements in finitely generated subgroups of $G$,
 or even in the group $G$ itself. In this case, as was mentioned in Modification 2,
 one can try to compute some   linear approximations of these lengths and then use
    heuristic algorithms to carry over LBA.

    In Section \ref{subsec:geodesic-braid} we discuss hardness of
    the problem of
    computing the geodesic length (GL problem)  in braid groups $B_n$ -- the
    original platforms of AAG protocol.  The time complexity of GLP in $B_n$
    relative to the standard set of Artin
    generators $\Sigma$
    is unknown.  We discuss some recent results and conjectures in this area.
    However, there are efficient linear  approximations
    of the geodesic length in $B_n$  relative to the set of generators $\Delta$
    (the generalized half-twists). Theoretically, this gives linear
    approximations of the geodesic length of elements in $B_n$ in the Artin
    generators, and, furthermore, linear approximations of geodesic
    inner length in quasi-isometrically embedded subgroups. If, as conjectured,  the
    set of quasi-isometrically embedded subgroups is exponentially
    generic in braid groups, then this gives a sound foundation for
    LBA in braid groups. Notice, that even  linear approximations
    alone are not entirely sufficient for successful  LBA. To get a
    precise solution of SCSP* one needs also a robust "local search" near
    a given  approximation of the solution. To this end several
    efficient  heuristic algorithms have been developed \cite{MSU1},
    \cite{MU3}. Nevertheless, by far none of them exploited directly the
    interesting interplay between geodesic lengths in $\Sigma$ and
    $\Delta$, as well as quasi-isometric embeddings of subgroups.

\subsection{Related algorithmic problems}
 \label{subsec:algor-problems-geodesic}

 We start with precise
    formulation of some problems related to computing geodesics in $G$.

\bigskip
\noindent
 {\bf Computing the geodesic length in a group (GL):}
Let $G$ be a group with a finite generating set $X$. Given an
element $w \in G$, as a product of generators form $X$,  compute the
geodesic length $l_X(w)$.

\bigskip
\noindent
 {\bf Computing the geodesic length in a subgroup (GLS):}
Let $G$ be a group with a finite generating set $X$ and $A$ a
subgroup of $G$ generated by a finite set of elements $Y = \{a_1,
\ldots,a_k\}$ of $G$ given as words from $F(X)$. Given an element $w
\in A$, as a product of generators of $A$,  compute the geodesic
length $l_Y(w)$.

\bigskip
There is another (harder) variation of this problem, that comes from
the SCSP* problem:

\bigskip
\noindent
  {\bf Computing the geodesic length in a subgroup (GLS*):}
Let $G$ be a group with a finite generating set $X$ and $A$ a
subgroup of $G$ generated by a finite set of elements $Y = \{a_1,
\ldots,a_k\}$ of $G$ given as words from $F(X)$. Given an element $w
\in A$, as a word from $F(X)$, compute the geodesic length $l_Y(w)$.

\bigskip
The following lemma is obvious. Recall, that The Membership
Search Problem (MSP) for a subgroup $A$ in $G$ requires for a given
element $w \in F(X)$, which belongs to $A$, to find a decomposition
of $w$ into a product of generators from $Y$ and their inverses.
\begin{lemma}
\label{le:MSP}
 Let $G$ be a finitely generated group and $A$ a
finitely generated subgroup of $G$. Then:
 \begin{itemize}
  \item [1)]  GLS is linear time reducible to GLS*;
   \item [2)] GLS* is linear time reducible to GLS relative to the
   Membership Search Problem in $A$.
 \end{itemize}
\end{lemma}

Observe, that if  GLS  has a "fast" solution for $A = G$ in $G$ then
there is a fast algorithm to find the geodesic length of elements of
$G$ with respect to $X$. In particular, the Word Problem in $G$ has
a fast decision algorithm. In some groups, like free groups or
partially commutative groups, given by the standard generating sets,
there are fast algorithms for computing the geodesic length of
elements. In many other groups, like braid groups, or nilpotent
groups, the computation of the geodesic length of elements is hard.
Nevertheless, in many applications, including cryptography,  it
suffices to have a fast algorithm to compute a reasonable, say
linear, approximation of the geodesic length of a given element.
 To this end we formulate the following problem.

\bigskip
\noindent
 {\bf Computing a linear approximation of the geodesic
length in a group (AGL):} Let $G$ be a group with a finite
generating set $X$. Given a word $w \in F(X)$ compute a linear
approximation  of the geodesic length of $w$. More precisely, find
an algorithm that for $w \in F(X)$ outputs a word $w' \in F(X)$ such
that $\lambda l_X(w) +c \geq l_X(w')$, where $\lambda$ and $c$ are
independent of $w$.

\bigskip
Another  problem is to compute a good approximation in a subgroup of
a group.

\bigskip
\noindent
 {\bf Computing a linear approximation of the geodesic
length in a subgroup  (AGLS):} Let $G$ be a group with a finite
generating set $X$  and $A$ a subgroup of $G$ generated by a finite
set of elements $Y = \{a_1, \ldots,a_k\}$ of $G$ given as words from
$F(X)$. Given an element $w \in A$, as a word from $F(X)$, compute
 a linear approximation of the geodesic length $l_Y(w)$ of $w$.

\bigskip

Assume now that there is a "fast"  algorithm to compute AGL in the
 group $G$. However, this does not imply that there is a fast
algorithm to compute a linear approximation of the geodesic length
in a given subgroup $A$ of $G$. Unless, the subgroup $A$ is
quasi-isometrically embedded in $G$.

\begin{lemma}
 Let $G$ be a group with a finite generating set $X$ and $\mathcal{A}$ is an algorithm to compute
 AGL in $G$ with respect to $X$. If $H$ is a quasi-isometrically
 embedded subgroup of $G$ generated by a finite set $Y$ then for every $w \in H$, given as a
 word from $F(X)$, the algorithm $\mathcal{A}$ outputs a word $w' \in F(X)$ such that $l_Y(w) \leq
 \mu l_X(w') + d$ for some constants $\mu$ and $d$ which depend only on $\mathcal{A}$ and $H$.
\end{lemma}
  \begin{proof} The proof is straightforward. \end{proof}

\subsection{Geodesic length in braid groups}
\label{subsec:geodesic-braid}

There is no any known efficient algorithm  to compute the geodesic
length of elements in braid groups with respect to the set $\Sigma$
of the standard Artin's generators. Some indications that this could
be a hard problem are given in \cite{PRaz}, where the authors prove
that the set of geodesics in $B_\infty$ is co-NP-complete. However,
in a given group, the problem of computing the length of a word
could be easier then the problem of finding a geodesic of the word.
Moreover, complexity of a set of geodesics in a group may not be a
good indicator of the time complexity of computing  the geodesic
length in a randomly chosen subgroup. In fact, it has been shown in
\cite{MSU1,MSU2} that in a braid group $B_n$ one can efficiently
compute a reasonable approximation of the length function on $B_n$
(relative to $\Sigma$) which gives a foundation for successful LBA,
without computing the length in the group.   Furthermore, there are
interesting open conjectures that, if settled affirmatively, will
lead to more efficient  algorithms for computing the length of
elements in braid groups and their subgroups. To explain this we
need to introduce some known facts and terminology.

The group $B_n$  has the classical Artin presentation:
\begin{displaymath}
B_n =
 \left\langle
\begin{array}{lcl}\sigma_1,\ldots,\sigma_{n-1} & \bigg{|} &
\begin{array}{ll}
 \sigma_i \sigma_j \sigma_i= \sigma_j \sigma_i \sigma_j & \textrm{if }|i-j|=1 \\
 \sigma_i \sigma_j=\sigma_j \sigma_i & \textrm{if }|i-j|>1
\end{array}
\end{array}
\right\rangle.
\end{displaymath}

By $l_\Sigma(w)$ we denote the length of a word $w \in B_n$ relative
to the generating set  $\Sigma = \{\sigma_1, \ldots,
\sigma_{n-1}\}$.

Elements in $B_n$ admit so-called {\em Garside} normal forms. These
forms are unique and the time complexity to compute the normal form
of an  element of $B_n$ given by a word $w\in F(\Sigma)$ is bounded
by $O(|w|^2 n^2)$. However, Garside normal forms are far from
being geodesic in $B_n$.

 In 1991 Patrick Dehornoy introduced in \cite{D1}
  the following notion of $\sigma$-positive braid word and a handle-reduction algorithm
  to compute a $\sigma$-positive representative of a given word. A braid word $w$
  is termed to be $\sigma_k$-positive (respectively,
negative), if it contains $\sigma_k$, but does not contain $\sigma_k
^{-1}$ and $\sigma_i ^{\pm 1}$  with $i < k$ (respectively, contains
$\sigma_k ^{-1}$, but not $\sigma_k$ and $\sigma_i ^{\pm 1}$ with $i
< k$).  A braid word $w$ is said to be $\sigma$-positive
(respectively, $\sigma$-negative), if it is $\sigma_k$-positive
(respectively, $\sigma_k$-negative) for some $k \leq  n-1$. A braid
word $w$ is said to be $\sigma$-consistent if it is either trivial
or $\sigma$-positive, or $\sigma$-negative.

\bigskip
{\bf Theorem.} [Dehornoy \cite{D1}]. {\it  For any braid $\beta \in
B_n$, exactly one of the following is true:

1) $\beta$ is trivial;

2) $\beta$ can be presented by $\sigma_k$-positive braid word for
some $k$;

3) $\beta$ can be presented by $\sigma_k$-negative braid word for
some $k$. }

In the latter two cases $k$ is unique.

\bigskip
 Thus, it makes sense to speak
about $\sigma$-positive and $\sigma_k$-positive (or $\sigma$-,
$\sigma_k$-negative) braids.

 The following question is of primary
interest when solving AGL in braid groups: is there a polynomial
$p(x)$ such that for every word $w \in F(\Sigma)$  $p(l_\Sigma(w))$
gives an upper bound  for the $\Sigma$-length of the shortest
$\sigma$-consistent braid word representing $w \in B_n$? Dehornoy's
original algorithms in \cite{D1}, and the handle reduction from
\cite{D2}),  and the algorithm from \cite{Fenn}, all of them  give
only an exponential bound on the length of the shortest
$\sigma$-consistent representative.

In \cite{Dyn} (see also \cite{D2,Fenn} for a  related  discussion)
Dynnikov and Wiest formulated the following

\begin{conjecture} There are numbers $\lambda , c $ such that
every braid $w \in B_n$ has a $\sigma$-consistent representative
whose $\Sigma$-length is bounded linearly by the $\Sigma$-length of
the braid.
 \end{conjecture}

They also showed that  the conjecture above has a  positive answer
 if  the $\Sigma$-length of elements is  replaced by the
$\Delta$-length (relative to a set of generators $\Delta$).

The set of generators $\Delta$ consists of the braids $\Delta_{ij},
1 \leq i < j \leq n,$ which are  the half-twists of strands $i$
through $j$:
 $$ \Delta_{ij} = (\sigma_i . . . \sigma_{j-1})(\sigma_i . . .
\sigma_{j-2}) . . . \sigma_i.$$
 $\Delta$ is a generating set of
$B_n$,  containing  the Artin's generators $\sigma_i =
\Delta_{i,i+1}$, and the Garside fundamental braid $\Delta_{1n}$.
 The {\it compressed} $\Delta$-length of a word $w$ of the form
 $$ w = \Delta_{i_1j_1} ^{k_1} . . .\Delta_{i_sj_s} ^{k_s},$$
 where $k_t \neq 0$ and $\Delta_{i_t,j_t} \neq  \Delta_{i_{t+1},j_{t+1}}$ for all $t$,
 is  defined by
 $$L_\Delta(w) = \Sigma_{i=1} \log_2(|k_i| + 1).$$
   For an element  $\beta \in B_n$  the value $L_\Delta(\beta)$ is
   defined by
    $$L_\Delta(\beta) = \min\{L_\Delta(w) \mid \  the \ word \ w \  represents \ \beta \}.$$
     Obviously, for any braid $\beta$, we have
      $$L_\Delta(\beta) \leq l_\Delta(\beta) \leq l_\Sigma(\beta).$$

The modified conjecture assumes the following extension of the
notion of $\sigma$-positive braid word: a word in the alphabet
$\Delta = \{\Delta_{ij} \mid 0<i<j<n \}$ is said to be
$\sigma$-positive if, for some $k < l$, it contains $\Delta_{kl}$,
and contains neither $\Delta_{ kj}^{-1}$ nor $\Delta_{ ij}^{\pm 1}$
with $i < k$ and any $j$. In other words, a word $w$ in letters
$\Delta_{ij}$ is $\sigma$-positive (negative) if the word in
standard generators $\sigma_i$ obtained from $w$ by the obvious
expansion is.

\bigskip
{\bf Theorem } [Dynnikov, Wiest \cite{Dyn}].  {\it Any braid $\beta
\in B_n$ can be presented by a $\sigma$-consistent word $w$ in the
alphabet $\{\Delta_{ij}\}$ such that $$ l_\Delta(w) \leq  30n
l_\Delta (\beta).$$ }

\medskip
This theorem gives a method  to approximate geodesic length in braid
groups, as well as in its quasi-isometrically embedded subgroups. It
remains to be seen whether this would lead to more efficient
versions of LBA or not.

\section{Quotient attacks}
\label{se:quotient-attacks}

In this section we describe a new type of attacks, which we term
{\em quotient attacks} (QA). In fact, the quotient attacks are just
fast generic algorithms to solve such search problems in groups as
the Membership Search Problem (MSP),  the Simultaneous Conjugacy
Search Problem (SCSP), the Simultaneous Conjugacy Search Problem
relative a to a subgroup (SCSP*), etc. The main idea behind QA is
that to solve a problem in a group $G$ it suffices, on most inputs,
to solve it in a quotient $G/N$, provided $G/N$ has generic free
basis property and a fast decision algorithm for the problem. In
particular, this is the case, if $G$ has a free non-abelian
quotient. Notice, that a similar idea was already exploited in
\cite{KMSS1}, but there the answer was given only for inputs in "No"
part of the decision problem, which, obviously, does not apply to
search problems. The strength of our approach comes from the extra
requirement that $G/N$ has the free basis property.

In Sections \ref{subsec:MP-free-groups} and
\ref{subsec:SCSP-free-groups} we discuss the Conjugacy and
Membership Problems  in all their variations  in free groups. Some
of these results were  known in folklore, some could be found in the
literature. Nevertheless, we sketch most of the proofs
here, since this will serve us as the base for solving similar
problems in other groups.

\subsection{Membership Problems  in free groups}
  \label{subsec:MP-free-groups}

In this section we discuss some algorithms to solve the  Membership
Problems in all their variations in free groups. We start with  the
classical Membership Problem (MP). Everywhere below $G$ is a fixed
group generated by a finite set $X$.

\bigskip
\noindent
 {\bf The Membership Problem (MP):} Let $A = \gp{a_1, \ldots, a_m}$ be
 a fixed finitely generated subgroup of $G$ given by a finite set of
 generators $a_1, \ldots, a_m$ (viewed as words in $F(X)$). Given a
 word $w \in F(X)$ decide whether  $w$ belongs to $A$ or not.

 \bigskip
 When the subgroup $A$ is not fixed, but comes as a part of the input (like in AAG scheme)
  then the problem is more precisely described in its {\em uniform}
  variation.

\bigskip
\noindent
 {\bf The Uniform Membership Problem (UMP):} Given a finite tuple of elements
 $w, a_1, \ldots, a_m \in F(X)$ decide whether or not $w$ (viewed as an element of $G$)
 belongs to the subgroup $A$ generated by
 the elements $a_1, \ldots, a_m$ in $G$.

 \bigskip
To solve MP in free groups we use the folding technique introduced
by Stallings in \cite{S}, see also \cite{KM} for a more detailed
treatment. Given a tuple of words
  $a_1, \ldots, a_m \in F(X)$ one can construct a finite deterministic automaton
  $\Gamma_A$, which accepts a reduced word $w \in F(X)$  if and only
  if $w$ belongs to the subgroup $A = \langle a_1, \ldots, a_m \rangle$ generated by
  $a_1, \ldots, a_m $ in $F(X)$.

To describe the time complexity of MP and UMP recall that for a
given positive integer $n$ the function $log_2^\ast n$ is defined as
the least natural number $m$ such that $m$-tower of exponents of $2$
exceeds $n$, or equivalently, $log_2 \circ log_2 \circ \ldots \circ log_2 (n)
\leq 1$, where on the left one has composition of $m$ logarithms.

\begin{lemma}\label{le:UMP-free}
There exists an algorithm which for any input
$w, a_1, \ldots, a_m \in F(X)$ for UMP finds the correct answer
in nearly linear time $O(|w|+ nlog^\ast n)$ where $n = \sum_{i=1}^k |a_i|$.
Furthermore, the algorithm works in linear time $O(|w|+ n)$
on exponentially generic set of inputs.
\end{lemma}

\begin{proof}
Indeed, given $w, a_1, \ldots, a_m \in F(X)$  one can construct
$\Gamma_A$ in worst time $O(nlog^\ast n)$ (see \cite{touikan}) and check if $\Gamma_A$
accepts $w$ or not in time $O(|w|)$, as required.

To prove the generic estimate recall that the set of $m$-tuples
$a_1, \ldots, a_m \in F(X)$ satisfying $1/4$-condition is exponentially generic
and the Stalling's procedure constructs the automaton $\Gamma_A$
in linear time $O(n)$.

\end{proof}

In cryptography, the search variations of MP and UMP  are the most
interesting.

\bigskip
\noindent
 {\bf The Membership Search Problem (MSP):}
 Let $A = \gp{a_1, \ldots, a_m}$ be
 a fixed finitely generated subgroup of $G$ given by a finite set of
 generators $a_1, \ldots, a_m$, viewed as words in $F(X)$. Given a
 word $w \in F(X)$, which belongs to $A$,  find a representation of $w$
 as a product of the generators
 $a_1, \ldots, a_m$ and their inverses.

\bigskip
\noindent
 {\bf The Uniform Membership Search Problem (UMSP):} Given a finite tuple of elements
 $w, a_1, \ldots, a_m \in F(X)$ such that $w \in A = \gp{a_1, \ldots, a_m}$ find
 a representation of $w$ as a product of the generators
 $a_1, \ldots, a_m$ and their inverses.

Time complexity upper bounds for MSP easily follow from the
corresponding bounds for $MP$.
\begin{lemma}
\label{le:MSP-free} The time complexity of MSP in a free group is
bounded from above by $O(|w|)$.
\end{lemma}
  \begin{proof}
  Let $A = \gp{a_1, \ldots, a_m}$  be a fixed finitely generated
  subgroup of $G$. As was mentioned above in time $O(nlog^\ast n)$, where
 $n = |a_1| + \ldots + |a_n|$, one can construct the Stallings'
 folding $\Gamma_A$. In linear time in $n$, using the breadth first
 search,  one can construct a Nielsen basis $S = \{b_1, \ldots,b_n\}$
 of $A$ (see \cite{KM}).
  Now, given a word $w \in F(X)$, that belongs to $A$, one can follow the accepting  path
  for $w$ in  $\Gamma_A$ and rewrite $w$ as a product of generators from $S$ and their inverses.
  This requires linear time in $|w|$. It is suffices to notice that
  the elements $b_i$ can be expressed as fixed products of
  elements from the initial generators of $A$, $b_i = u_i(a_1,
  \ldots,a_n)$, $i = 1, \ldots, m$, therefore any expression of $w$
  as a product of elements from $S^{\pm 1}$ can be rewritten in a
   linear time into a product of the initial generators.
  \end{proof}

Observe, that in the proof above we used the fact that  any  product
of new generators $b_i$ and their inversions can be rewritten in
linear time into a product of the old generators $a_i$ and their
inversions. That held because we assumed that one can rewrite the
new generators $b_i$ as products of the old generators $a_i$ in a
constant time. This is correct if the subgroup $A$ is fixed.
Otherwise, say in UMSP,  the assumption does not hold anymore. It is
not even clear whether one  can do it in polynomial time or not. In
fact, the time complexity of UMSP is unknown.  The following problem
is of prime interest in this area.

\begin{problem}
 \label{pr:UMSP}
 Is the time complexity of UMSP in  free groups polynomial?
\end{problem}

However, the generic case complexity of UMSP in free groups is
known.

\begin{lemma} \label{le:UMSP-free-generic}
The generic case time complexity of UMSP in free groups is linear.
More precisely, there is an exponentially
generic subset $T \subseteq F(X)^n$ such that for every tuple $(w,
a_1, \ldots, a_m) \in F(X) \times T$, such that $w \in \langle a_1,
\ldots, a_m\rangle$, one can express $w$ as a product of $a_1,
\ldots, a_m$ and their inverses in time $O(|w|+n)$ where
$n = |a_1| + \ldots + |a_n|$.
\end{lemma}
  \begin{proof}
 Notice, first, that if in the argument of Lemma \ref{le:MSP-free}
 the initial set of generators $ a_1, \ldots, a_m$ of a subgroup $A$
 satisfy $1/4$-condition then the set of the new generators $b_1, \ldots,
 b_m$ coincides with the set of the initial generators (see
 \cite{KM} for details). Moreover, as was noticed in the proof of
 Theorem \ref{th:LBA-free-groups}
  the set $T$ of tuples $(a_1, \ldots, a_m) \in F(X)^m$,
  satisfying $1/4$-condition
   is exponentially generic. Hence the argument from  Lemma
   \ref{le:MSP-free}  proves the required upper bound for UMSP on
   $T$.
 \end{proof}

\subsection{The Conjugacy  Problems  in free groups}
  \label{subsec:SCSP-free-groups}

Now we turn to the conjugacy problems in free groups. Again,
everywhere below $G$ is a fixed group generated by a finite set $X$.

It is easy to see that the CP and CSP in free groups are decidable
in at most quadratic time. It is quite tricky to show that CP and
CSP are decidable in free groups in linear time! This result
is based on Knuth-Morris-Pratt substring searching algorithm \cite{KMP}.
Similarly, the Root Search Problem (listed below) is
decidable in free groups in linear time.

\bigskip
\noindent
 {\bf The Root Search Problem (RSP):} Given a word  $w \in F(X)$
 find a shortest word $u \in F(X)$ such that $w = u^n$ for some
 positive integer $n$.

 \bigskip
Notice, that  RSP in free groups can be interpreted as a problem of
finding a single  generator of the centralizer  of a  non-trivial
element.

\begin{theorem} \label{le:SSCP-free-groups}
The Simultaneous Conjugacy Problem (SCP)  and Simultaneous Conjugacy
Search Problem  (SCSP)   are in linear time reducible to CP, CSP,
and RP in free groups. In particular, it is decidable in linear time.
\end{theorem}

\begin{proof}
 We briefly outline an  algorithm that simultaneously solves the problems
 SCP and SCSP in free groups, i.e.,
 given a  finite system of conjugacy equations
 \begin{equation}
 \label{eq:system}
 \left\{
 \begin{array}{l}
    u_1^x = v_1, \\
    \ldots \\
    u_n^x = v_n, \\
 \end{array}
 \right.
 \end{equation}
the algorithm decides whether or not this system has a solution in a
free group $F(X)$, and if so, it finds a solution. Using the
decision algorithm for CP one can check whether or not there is an
equation in (\ref{eq:system}) that does not have solutions in $F$.
If so the whole system does not have solutions in $F$ and we are
done.
  Otherwise, using the algorithm to solve CSP in $F$ one can find a
  particular solution $d_i$ of every equation $ u_i^x = v_i$ in (\ref{eq:system}).
  In this case the set of all solutions of the equation $ u_i^x = v_i$
  is equal to the coset $C(u_i)d_i$ of   the centralizer $C(u_i)$.
  Observe, that using the decision algorithm for RSP one can find
  a generator (the root of $u_i$) of the
  centralizer $C(u_i)$ in $F$.

Consider now the first two equations in (\ref{eq:system}).   The
system
  \begin{equation}
  \label{eq:system-abridged}
  u_1^x = v_1, u_2^x = v_2
  \end{equation}
  has a solution in $F(X)$ if and only if the
intersection $V = C(u_1)d_1 \cap C(u_2)d_2$ is non-empty. In this
case
 $$V = C(u_1)d_1 \cap C(u_2)d_2 = \left ( C(u_1) \cap C(u_2) \right )d$$
 for some $d \in F$.

 If $[u_1,u_2] = 1$ then $V$, as the intersection of two cosets, is
 non-trivial if and only if the cosets coincide, i.e.,
 $[u_1,d_1 d_2^{-1}] = 1$. This can be checked in linear time
 (since the word problem in $F(X)$ is in linear
 time). Therefore, in linear time we either check that the system
 (\ref{eq:system-abridged}), hence the system (\ref {eq:system}),
  does not have solutions at all, or we confirm that (\ref{eq:system-abridged})  is equivalent to one
 of the equations, so (\ref {eq:system}) is equivalent to its own subsystem,
 where the first equation is
 removed. In the latter case induction finishes the proof.

If $[u_1,u_2] \neq 1$ then $C(u_1) \cap C(u_2) = 1$,  so either $V =
\emptyset$ or $V = \{d\}$,  in both cases one can easily find all
solutions of (\ref{eq:system}). Indeed, if $V = \emptyset$ then
(\ref{eq:system})  does not have solutions at all. If $V = \{d\}$,
then $d$ is the only potential solution of (\ref{eq:system}), and
one can check whether or not  $d$ satisfies all other equations in
(\ref{eq:system}) in linear time  by the direct verification.

 Now the problem is  to verify in linear time whether $V = \emptyset$ or
 not, which is equivalent to solving an equation
  \begin{equation}
  \label{eq:powers-equation}
  u_1^md_1 = u_2^kd_2
   \end{equation}
 for integers $m,k$. Finding in linear time the cyclically reduced
 decompositions of $u_1$ and $u_2$ one can rewrite the equation
 (\ref{eq:powers-equation}) into an equivalent one in the form:
 \begin{equation}
  \label{eq:powers-equation-two}
   w_2^{-k} c w_1^m = b
  \end{equation}
where $w_1, w_2$ are cyclically reduced forms of $u_1, u_2$, and
either $w_2^{-1}c$ or $cw_1$ (or both) are reduced as written, and
$b$ does not begin  with $w_2^{-1}$ and does not end with $w_1$.
Again, in linear time one can find the maximal possible cancelation
in $w_2^{-k} c$, and in  $cw_1$, and rewrite
(\ref{eq:powers-equation-two}) in the form:
\begin{equation}
  \label{eq:powers-equation-three}
   w_2^{-k}  \tilde{w}_1^s = \tilde{b}
  \end{equation}
where $\tilde{w}_1$ is a cyclic permutation of $w_1$, and
$|\tilde{b}| \leq |b| + |w_1|$. Notice, that
 two cyclically reduced periodic words $w_2, \tilde{w}_1$ either commute or do not
have a common subword of length exceeding $|w_2| + |\tilde{w}_1|$.
If they commute then the equation (\ref{eq:powers-equation-three})
becomes  a power equation, which is easy to solve. Otherwise,
executing (in linear time) possible cancelation  in the left-hand
side of (\ref{eq:powers-equation-three}) one arrives to an equation
of the type
\begin{equation}
   w_2^{-r}  e \tilde{w}_1^t = \tilde{b}
  \end{equation}
where there is no cancelation at all. This can be easily solved for
$r$ and $t$.  This proves the result.

\end{proof}

As we have seen in the proof of Theorem \ref{le:SSCP-free-groups}
one of the main difficulties in solving SCSP  in groups  lies in
computing the intersection of two finitely generated subgroups or
their cosets. Notice, that finitely generated subgroups of $F(X)$
are  regular sets (which are accepted by their Stallings' automata).
It is well known in the language theory that the intersection of two
regular sets is again regular, and one can find an automaton
accepting the intersection in at most quadratic time. This leads to
 the following corollary.

\begin{corollary} \label{co:SSCP-free-groups}
The SCSP* in free groups is decidable in at most quadratic time.
\end{corollary}

\begin{proof}
Recall from the proof of Theorem \ref{le:SSCP-free-groups} that the
algorithm solving  a finite system of conjugacy equations in a free
group either decides that there is no solution to the system, or
produces a unique solution, or gives the whole solution set as a
coset $Cd$ of some  centralizer $C$.  In the first case,  the
corresponding SCSP* has no solutions in a given finitely generated
subgroup $A$; in the second case, given a unique solution $w$ of the
system one can construct the automaton $\Gamma_A$, that accepts $A$,
and check whether $w$ is in $A$ or not (it requires $nlog^\ast n$
time); and in the third case, one needs to verify if $Cd \cap A$ is
empty or not - this can be done, as we have mentioned above,  in at
most quadratic time (as the intersection of two regular subsets).
\end{proof}

Observe from the proof above, that the most time consuming case in
solving SCSP* in free groups occurs when all the elements $u_1,
\ldots, u_n$ in the system (\ref{eq:system}) commute. The set of
such inputs for SCSP* is, obviously,  exponentially negligible.
As we proved in Theorem \ref{th:LBA-free-groups} that LBA relative to $l_T$
solves SCSP* in linear time.

Since AAG  is reducible in linear time to SCSP* (Lemma
\ref{le:AAG-reduced}) we have the following results.
\begin{corollary}
\label{co:AAG-free} The following hold in an arbitrary free group
$F$.

  \begin{itemize}
  \item [1)]  The AAG algorithmic problem in $F$ is decidable in at most
  quadratic time in the size of the input (the size of the public
  information in the AAG scheme).
   \item [2)] The AAG algorithmic problem in $F$ is decidable in
   linear time on an exponentially generic set of inputs.
   \end{itemize}
\end{corollary}

\subsection{The MSP and SCSP* problems in groups with "good" quotients}
\label{subsec:MSP-SCSP-quotients}

In this section we discuss the generic complexity of the Membership
Search Problem MSP and the Simultaneous Conjugacy Search Problem
relative to a subgroup SCSP* in groups that have "good" factors in
$\mathcal{FB}_{exp}$.

Let $G$ be a group generate by a finite set $X$, $G/N$ is a quotient
of $G$, and $\phi:G \rightarrow G/N$ a canonical epimorphism. Let $H
= \gp{u_1, \ldots, u_k}$ be a finitely generated subgroup of $G$. To
solve the membership search problem for $H$ one can employ the
following simple heuristic idea which we formulate as an algorithm.

\begin{algorithm}\label{al:Heur_MSP}{\bf(Heuristic solution to MSP)}
    \\{\sc Input:}
A word $w = w(X)$ and generators $\{u_1, \ldots, u_k\} \subset F(X)$ of a subgroup $H$.
    \\{\sc Output:}
A representation $W(u_1, \ldots, u_k)$ of $w$ as an element of $H$
or $Failure$.
    \\{\sc Computations:}
\begin{itemize}
    \item[A.]
Compute the generators $u_1^\phi, \ldots, u_k^\phi$ of $H^\phi$ in
$G/N$.
    \item[B.]
Compute $w^\phi$, solve MSP for $w^\phi$ and $H^\phi$,  and find a
representation $W(u_1^\phi, \ldots, u_k^\phi)$ of $w^\phi$ as a
product of the generators of $u_1^\phi, \ldots, u_k^\phi$ and their
inverses.
    \item[C.]
Check if $W(u_1, \ldots, u_k)$ is equal to $w$ in $G$. If this  is
the case then output $W$. Otherwise output $Failure$.
\end{itemize}
\end{algorithm}

Observe that  to run Algorithm \ref{al:Heur_MSP} one needs to be
able to solve MSP in the quotient $G/N$ (Step B) and to check the
result in the original group (Step C), i.e., to solve the Word
Problem in $G$. If these conditions are satisfied Algorithm
\ref{al:Heur_MSP} is a partial  deterministic correct algorithm, it
gives only the correct answers. However, it is far from being
obvious, even the conditions are satisfied,  that this heuristic
algorithm can be robust in any interesting class of groups. The next
theorem, which is the main result of this section,  states that
Algorithm \ref{al:Heur_MSP} is very robust  for groups from
$\mathcal{FB}_{exp}$ with  a few additional requirements.

\begin{theorem}
\label{th:generic-MSP}
 {\bf (Reduction to a quotient)}
Let $G$ be a  group generated by a finite set $X$  and with the Word
Problem in a complexity class $C_1(n)$. Suppose $G/N$ is a quotient
of $G$ such that:
 \begin{itemize}
    \item [1)]
$G/N \in \mathcal{FB}_{exp}$.
    \item [2)]
The canonical epimorphism $\phi:G \rightarrow G/N$ is computable
within time $C_2(n)$.
    \item [3)]
For every $k \in \mathbb{N}$ there exists an algorithm
$\mathcal{A}_k$ in a complexity class $C_3(n)$,
which solves the Membership Search  Problem in $G/N$ for
an exponentially generic set $M_k \subseteq F(X)^k$ of descriptions of
$k$-generated subgroups in $G/N$.
\end{itemize}
  Then for every $k$ Algorithm \ref{al:Heur_MSP} solves
  the Membership Search Problem on an exponentially generic set $T_k \subseteq F(X)^k$
  of descriptions of $k$-generated subgroups in $G$. Furthermore,  Algorithm \ref{al:Heur_MSP}
  belongs to the complexity class $C_1(n)+C_2(n)+C_3(n)$.
\end{theorem}
  \begin{proof}
We need to show that Algorithm \ref{al:Heur_MSP} successfully
halts on an exponentially generic set of tuples from $F(X)^k$.
 By the conditions of the
  theorem the set $S_k$ of all $k$-tuples from
  $F(X)^k$ whose images in $G/N$ freely generate free subgroups is
  exponentially generic, as well as, the set $M_k$ of all tuples
  from $F(X)^k$ where the algorithm $\mathcal{A}_k$ applies. Hence
  the intersection $T_k = S_k \cap M_k$ is exponentially generic in
  $F(X)^k$.  We claim that Algorithm \ref{al:Heur_MSP} applies to
  the subgroups with descriptions from $T_k$. Indeed, the algorithm
  $\mathcal{A}_k$ applies to subgroups generated by tuples $Y = (u_1, \ldots,u_k)$
  from $T_k$, so if $w^\phi  \in H^\phi = \langle Y^\phi \rangle$
  then $\mathcal{A}_k$ outputs a required representation $w^\phi =
  W(Y^\phi)$ in $G/N$. Notice, that  $H^\phi$ is freely generated by
  $Y^\phi$ since $Y \in S_k$, therefore $\phi$ is injective on $H$.
  It follows that $w =W(Y)$ in $G$, as required. This proves the
  theorem.
  \end{proof}

 Theorems \ref{th:LBA-FB} and  \ref{th:generic-MSP}  imply the
 following result.

\begin{corollary}\label{co:co:generic-LBA}
Let $G$ be as in Theorem \ref{th:generic-MSP}. Then for every $k,m
>0$ there exists an algorithm $\mathcal{C}_{k,m}$ that solves the
SCSP* on an exponentially generic subset  of the set of all inputs
$I_{k,m}$ for SCSP*. Furthermore, $\mathcal{C}_{k,m}$ belongs to
the complexity class $n^2 + C_1(n) + C_2(n) + C_3(n)$.
\end{corollary}

\begin{corollary}\label{co:generic-LBA-Pure-braids}
Let $G$ be a group of pure braids $PB_n$, $n\geq 3$,  or a
non-abelian partially commutative group $G(\Gamma)$.  Then for
every $k,m >0$ there exists an algorithm $\mathcal{C}_{k,m}$ that
solves the SCSP* on an exponentially generic subset of the set of
all inputs $I_{k,m}$ for SCSP*. Furthermore, $\mathcal{C}_{k,m}$
belongs to the complexity class $O(n^2)$.
\end{corollary}

\begin{proof}
Recall that the for any pure braid group or a non-abelian partially
commutative group the Word problem can be solved by a quadratic time
algorithm.  Now the statement follows from Corollary
\ref{co:co:generic-LBA} and Corollaries \ref{co:pure-braids} and
\ref{co:part-comm}.

\end{proof}

\end{document}